\documentclass[12pt, leqno, twoside]{amsart}
\usepackage{amsmath, amsthm, amscd, amsfonts, amssymb, graphicx, color}
\usepackage{mathrsfs}
\usepackage[bookmarksnumbered, colorlinks, plainpages]{hyperref}
\hypersetup{colorlinks=true,linkcolor=red, anchorcolor=green, citecolor=cyan, urlcolor=red, filecolor=magenta, pdftoolbar=true}
\usepackage[pagewise]{lineno}

\usepackage[]{draftwatermark}
\SetWatermarkText{DRAFT}
\SetWatermarkScale{6}

\newcommand{\disp}{\displaystyle} 
\usepackage{enumerate}
\newtheorem{theorem}{Theorem}[section]
\newtheorem{lemma}[theorem]{Lemma}
\newtheorem{proposition}[theorem]{Proposition}
\newtheorem{corollary}[theorem]{Corollary}
\theoremstyle{definition}

\theoremstyle{remark}
\newtheorem{remark}[theorem]{Remark}
\numberwithin{equation}{section}

\begin{document}

\setcounter{page}{1}

\title[Parabolic frequency on conformal Ricci flow]{Parabolic frequency monotonicity on the conformal Ricci flow}
 
\author{Abimbola Abolarinwa}
\address{Department of Mathematics, University
of Lagos, Akoka, Lagos State, Nigeria}
\email{a.abolarinwa1@gmail.com}


\author{Shahroud Azami}
\address{Department of  Pure Mathematics, Faculty of Science, Imam Khomeini International University,
Qazvin, Iran. }
\email{azami@sci.ikiu.ac.ir}

\subjclass[2010]{53C21,  53E20, 35K65,58J35}
\keywords{Frequency functional; Conformal Ricci flow; Drifting Laplacian;  Monotonicity;weighted measure}
 \date{April 6, 2023}

\begin{abstract}
This paper is devoted to the investigation of the monotonicity of parabolic frequency functional under conformal Ricci flow defined on a closed Riemannian manifold of constant scalar curvature and dimension not less than $3$.  Parabolic frequency functional for solutions of certain linear heat equation coupled with conformal pressure  is defined and its  monotonicity under the conformal Ricci flow is proved by applying  Bakry-\'Emery Ricci curvature bounds. Some consequences of the monotonicity are also presented.
\end{abstract}

\maketitle


 \section{Introduction and main results}

Frequency  functional (elliptic and parabolic) monotonicity has been a hot topic in the analysis of partial differential equations and differential geometry since the classical paper \cite{Am} owing to its several applications. On the other hand, geometric flows of time dependent Riemmanian metrics such as the Ricci flow and its various extensions and modification have been widely studied since the seminal paper \cite{Ha82} for their numerous topological, geometric and physical applications.  For better and improved topological and geometric results,  conformal Ricci flow has been introduced in a recent work \cite{Fi} as a  modification of the Ricci flow.  The principal aim of the present paper is to extend monotonicity of parabolic frequency functional  for linear heat equation to the setting of compact (without boundary) Riemannian maniofld evolving by the conformal Ricci flow,  and then investigate what geometric condition(s) is/are required for such monotonicity and as well as its possible applications. The precise definitions, history, applications of, and relevant literature on frequency functional and conformal Ricci flow are discussed in what follows.

\subsection{Frequency functionals}
Let $q$ be a fixed  point in $\mathbb{R}^{n}$ and $\mathcal{B}(q,r)$ be a ball  of radius $r>0$ centred at $q$.  Almgren \cite{Am} introduced the frequency functional  (which is known in literature as elliptic frequency functional) for a harmonic function $v(x)$ on $\mathbb{R}^{n}$, i.e., $\Delta_{\mathbb{R}^n} v(x)=0 \ \ \text{on}\ \mathbb{R}^n$, as follows
\begin{align}\label{F1}
E_F(r)=\frac{r\int_{\mathcal{B}(q,r)}|\nabla v(x)|^{2}dx}{\int_{\partial \mathcal{B}(q,r)}|v(x)|^{2}dS},
\end{align}
where $\partial \mathcal{B}(q,r)$ and $dS$ are respectively the boundary of $\mathcal{B}(q,r)$  and the induced $(n-1)$-dimensional Hausdorff measure on $\partial \mathcal{B}(q,r)$.  Here, $\nabla$ and $\Delta_{\mathbb{R}^n}$ are gradient and Laplace operators on $\mathbb{R}^n$, respectively.
The rate  of growth of harmonic function $v$  near the fixed point $q$ is determined by the functional in \eqref{F1}.  Furthermore,  Almgren \cite{Am} proved that $E_F(r)$ is monotone nondecreasing  with respect to $r$,  the consequence of which led to the study of the local regularity of harmonic functions and  minimal surfaces.  

Since the work of Almgren \cite{Am},  the monotonicity of elliptic functional $E_F(r)$ has been successfully applied in the analysis of more general elliptic and parabolic partial differential equations, and there have been also considerable generalization  to Riemannian manifolds.  We mention but a few literature: Garofalo and Lin \cite{GL1, GL2} investigated  the unique continuation properties  for elliptic operators by using the monotonicity of frequency functionals on Riemannian manifolds.  The authors in \cite{HL,Li,Lo} applied monotonicity of the frequency functionals to estimate the size of nodal and critical sets  of solutions  to elliptic and parabolic equations.  Colding and Minicozzi in \cite{CM1} applied  frequency monotonicity to prove finite dimensionality of the space of polynomial growth of harmonic functions on manifolds with nonnegative Ricci curvature and Euclidean volume growth, while in \cite{CM2} they extended the result to the case of static manifold using drifting Laplacian.
The counterpart of $E_F(r)$  for solutions to the heat equation  on $\mathbb{R}^{n}$  is called parabolic frequency  functional,  which was  first introduced by Poon \cite{P} to  the  study of  the unique  continuation   of solutions  to  parabolic equations on $\mathbb{R}^n$.   Consider a smooth solution $u=u(x,t)$ to the heat equation
\begin{align}\label{h}
\partial_t u -\Delta_M u=0 \ \ \ \text{in} \ M \times [0,T],
\end{align}
where $(M,g)$ is a complete  Riemannian manifold, $\Delta_M$ is the Laplace-Beltrami operator on $M$ and $T>0$. The Parabolic frequency for $u$ is defined by \cite{N} (see also \cite{P}) 
 \begin{align}
 P_F(t) = t\cdot \frac{\int_M |\nabla u|^2(x,T-t) H(x,y;t)d\mu(g)}{\int_Mu^2(x,T-t) H(x,y;t)d\mu(g)},
 \end{align}
where $H(x,y;t)$ is the fundamental solution to the heat equation \eqref{h}, $y$ being a referenced point  (not so important) in $M$,  and $d\mu(g)$ is the volume form with respect to the Riemannian metric $g$. Restricting $M$ to possessing nonnegative sectional  (or bisectional for holomorphic function) curvature and parallel Ricci curvature,   Poon \cite{P} and Ni \cite{N} proved that $P_F(t)$ is monotone nondecreasing by using Hamilton's matrix Harnack estimate \cite{Ha93}.  Recently, these results have been generalized to more general Riemannian manifolds by Li and Wang \cite{ XL}.

Let  $\tau(t)$ be the backward time, $\kappa(t)$ be the time-dependent function and $d\nu$ be the
weighted measure.  For  a solution $u(t)$
of the heat equation,  Baldauf and Kim \cite{BJ} defined the  parabolic frequency  as follows
$$U(t)=-\frac{\tau(t)||\nabla u||_{L^{2}(d\nu)}^{2}}{|| u||_{L^{2}(d\nu)}^{2}}e^{-\int \frac{1-\kappa(t)}{\tau(t)}dt}.$$
In the above definition, the exponential  term involving time-dependent function $\kappa$ serves as a correction term which depends on the geometry of the flow, analogous to the error term involving $r$ in the elliptic case. The authors \cite{BJ} proved that parabolic frequency $U(t)$ for the solution of heat equation is monotone increasing along the Ricci flow with bounded Bakry-\'{E}mery Ricci curvature.  See also the recent preprints \cite{AA,AA1,BL,LLX} for related results under the Ricci-Bourguignon,  Ricci-harmonic and  mean curvature flows.  Motivated by the above cited works, we study monotonicity of a well defined parabolic frequency function (see \eqref{Q} below) for a form of linear heat equation defined in \eqref{lh} along the conformal Ricci flow.   This study is more interesting since conformal Ricci flow performs better than (and even complementary to) Ricci flow in searching for certain geometric features,  and has wider applications in conformal geometry of constant scalar curvature.

\subsection{Conformal Ricci flow}
The conformally modified Ricci flow was introduced by Fischer in \cite{Fi} and named conformal Ricci flow as a result of the role played by conformal geometry in the derivation of its equations.   Precisely, let $(M,g_0)$ be a smooth $n$-dimensional $(n\ge 3)$ closed connected  manifold  together with Riemannian metric $g_0$ of constant scalar curvature $R_0$. The conformal Ricci flow is defined by a one-parameter family of metric $g(t)$ satisfying the following parabolic system
\begin{equation}\label{CRF}
\begin{cases}
\displaystyle  \frac{\partial g(t)}{\partial t}+2\left(Ric(t)-\frac{R_0}{n}g(t) \right) = -2p(t)g(t) , \ &(x,t) \in M\times (0,T),  
\\ 
\quad  R_{g(t)}=R_0, \ & (x,t) \in M\times [0,T),
\end{cases}
\end{equation}
together with the initial condition $g(0)=g_0$ and a family of function $p(t)$, $t\in [0,T)$, where $Ric(t)$ and $R_{g(t)}$ are the Ricci tensor and scalar curvature of the evolving metric $g(t)$, respectively. 
By the constraint equation $R_{g(t)}=R_0$ in system \eqref{CRF}, the flow is known to preserve constant scalar curvature of the evolving metric. Indeed, this accounts for naming the function, $p=p(t)$, conformal pressure, since it serves as time-dependent Lagrange multiplier and makes the term $-p(t)g(t)$ acting as the constraint force necessary to preserve the scalar curvature constraint. Consequently, $p(t)$ is known to solve a time-dependent elliptic partial differential equation under the flow
\begin{equation}\label{e12}
(n-1)\Delta p + R_0p = -\left|Ric-\frac{R_0}{n}g\right|^2   \hspace{1cm} \text{in} \ \ M\times [0,T).
\end{equation}
Considering the role of the conformal pressure, the function $p(t)$ is expected to be zero at an equilibrium point and strictly positive otherwise. Hence, the equilibrium points of the conformal Ricci flow are characterized by Einstein metrics, and  the term $-2(Ric-\frac{R_0}{n}g)$ can then be regarded as a measure of deviation of the flow from an equilibrium point. Conformal Ricci flow as a quasi-linear parabolic-elliptic system (from view point of \eqref{CRF} and \eqref{e12}) is analogous to semi-linear  Navier-Stokes equations for incompressible viscous fluid flow.

Since the volume of a Riemannian manifold $(M,g)$ is a positive real number and the scalar curvature is a real-valued function on $M$, the constraint on $R_{g(t)}$ is considerably more drastic than the volume constraint of the Ricci flow \cite{Ha82,Ha93}.  The Ricci flow in general does not preserve the property of constant scalar curvature. Thus, the configuration space of the conformal Ricci flow equations is considerably smaller than that of the Ricci flow.  From geometric point of view,  working on a smaller configuration space may be more advantageous than working on a larger configuration space (see further discussion in \cite[Section 1.6]{Fi}).  More concrete similarities and differences between conformal Ricci flow and the classical Ricci flow, as well as some possible applications of conformal Ricci flow to $3$-manifold geometry and conformal geometry are highlighted in \cite{Fi}.  Fischer's paper \cite{Fi} has also presented a proof of the short-time existence and uniqueness of the conformal Ricci flow on closed manifolds with negative constant scalar curvature $R_0<0$.  In that same paper, he also observed that Yamabe constant is strictly increasing along the flow on negative Yamabe type closed manifolds. For detail discussion on Yamabe problem, interested readers can consult the book by Aubin \cite{Au}. The following references \cite{BB,Bel,LASA, LLWY,LQZ, LQZ1,SZ} can be found for further studies on conformal Ricci flow.

\subsection{Main results}
Denote the partial derivative of any time-dependent quantity by $\partial_t$ (i.e., $\partial_tu(x,t) = u_t(x,t)$). Consider a one-parameter family of metrics $g = g(t)$,  $t \in [0, T )$, $T>0$, on an $(m+1)$-dimensional closed manifold with initial metric $g(0)$ having constant scalar curvature $-m(m+1)$ which is preserved under the flow as  $R_{g(t)}\equiv -m(m+1)$.  Referring to \eqref{CRF} and \eqref{e12}, one sees that $(M^{m+1},g(t),p(t))$ evolves by conformal Ricci flow given in  the following system
\begin{align}\label{CRF2}
\left\{
\begin{array}{ll}
\displaystyle \partial_t g = -2(Ric(t) + (m+p(t))g(t))  \, \, \, \,  & \text{on}\ \ M^{m+1}\times [0, T), \\ 
\displaystyle (-\Delta_{g(t)} + (m+1))p(t) = \frac{1}{m} |Ric +mg|^2 \, \, \, \,  & \text{on} \ \  M^{m+1}\times [0, T),  \\
g(0)  = g_0.
\end{array}
\right.
\end{align}
Denote the Laplace-Betrami operator with respect to $g(t)$ by $\Delta_{g(t)}$.  For a smooth function $v$ and a time-dependent nonnegative scalar non-dynamical field (conformal pressure) $p(t)$ satisfying the second equation in system \eqref{CRF2}, we consider the following linear heat equation
\begin{equation}\label{lh}
v_{t}(x,t)-\Delta_{g(t)} v(x,t)=\bar{p}(t)v(x,t),  \ \ \ v(0,x)=v_0,
\end{equation}
where $\disp \bar{p}(t)=\max_{x\in M}p(t)$, on $M\times [0,T)$.  Note that the conformal pressure $p(t)$ is finite \cite{Fi,SZ} (see also Lemma \ref{lemp} below), and $\bar{p}(t)$ can be thought of as its upper bound. Thus, $\bar{p}(t)$ is viewed as material-constant (or space-constant) pressure function.  Following \cite{LLX} we define the parabolic frequency functional $Q(t)$ for the  solutions of heat equation  \eqref{lh} along the  conformal Ricci flow  as follows
 \begin{equation}\label{Q}
Q(t)=\frac{h(t)\int_{M}|\nabla_{g(t)} v|_{g(t)}^{2}dV_{g(t)}}{\int_{M}v^{2}dV_{g(t)}}e^{-\int_{t_{0}}^{t} \left(2\bar{p}(s)+\frac{h'(s)+k(s)}{h(s)}\right)ds}\
\end{equation}
where $h$ and $k$ are  smooth functions with respect to time-variable $t \in [t_{0},t]\subset[0,T)$, and $dV_{g(t)}$  is the weighted measure on $(M,g(t))$ (see the appropriate definition of $dV_{g(t)}$ in (\eqref{24}) below).   The involvement of a finite time-dependent scalar function $\bar{p}$ in the exponential term of the functional is natural: First, to reflect the coupling of the solution of the heat equation with the conformal pressure, and secondly, to show the conformal nature of $Q(t)$ in preserving the scalar curvature constraint. (See \cite{Fi} for complete description of the pressure field $p(t)$).

In this paper we prove monotonicity of \eqref{Q} along the flow \eqref{CRF2} coupled with the heat equation \eqref{lh} and then obtain our main result.  Let $\mathscr{R}ic_f$ and $\mathscr{L}_f$ be Bakry-\'Emery Ricci curvature tensor and drifting Laplacian, respectively.
 
\begin{theorem}\label{thm11}
Let $(M^{m+1},g(t),p(t)),\,\,t\in[0,T)$ be  a solution to the conformal Ricci flow \eqref{CRF2} with $\mathscr{R}ic_f\leq (\frac{k(t)}{2h(t)}+\frac{R_g}{m+1})g$ and $R_g=-m(m+1)$. 
\begin{enumerate}
\item [(i)] If $h(t)$ is positive (i.e., $h(t)>0$), then the parabolic frequency $Q(t)$  is monotone nonincreasing along the conformal Ricci flow.
\item[(ii)] If $h(t)$ is  negative (i.e., $h(t)<0$), then the parabolic frequency $Q(t)$  is monotone nondecreasing along the  conformal Ricci flow.
\end{enumerate}
Furthermore,  $Q'(t)=0$ only if $v$ is an eigenfunction of $\mathscr{L}_f$ satisfying $-\mathscr{L}_f v=c(t)v$. 
\end{theorem}
As an application we have the following corollaries.
\begin{corollary}\label{cor2}(Backward uniqueness).
 Assuming the hypotheses of Theorem \ref{thm11} hold. If $h(t)<0$ and $v(\cdot,b)=0$, then $v(\cdot,t)=0$ for any $t\in [a,b) \subset (0,T)$, $a<b$.
\end{corollary}

\begin{corollary}\label{cor3} (Eigenvalue monotonicity).
Define the first nonzero eigenvalue of $(M^{m+1},g(t),p(t))$ with respect to the drifting Laplacian on weighted measure $dV_{g(t)}$ by
\begin{align*}
\lambda(g(t)):= \inf \frac{\int_{M} -\langle \mathscr{L}_f u, u\rangle dV_{g(t)}}{\int_{M} u^2 dV_{g(t)}}:  \ u \in W^{1,2}(M^{m+1})\setminus\{0\} \  \int_M u dV_{g(t)}=0,
\end{align*}
where $u(t)$ solves the linear heat equation \eqref{lh}.
Suppose $(M^{m+1},g(t),p(t)),\,\,t\in[0,T)$ solves  \eqref{CRF2} with $\mathscr{R}ic_f\leq (\frac{k(t)}{2h(t)}+\frac{R_g}{m+1})g$ and $R_g=-m(m+1)$.   Then for any $t\in [t_0,t_1]\subset (0,T)$
\begin{enumerate}
\item [(i)] If $h(t)>0$ then $h(t)\lambda(g(t))$  is a monotone decreasing function.
\item[(ii)] If $h(t)<0$ then $h(t)\lambda(g(t))$  is a monotone increasing function. 
\end{enumerate}
\end{corollary}

\begin{remark}
The Bakry-\'Emery Ricci condition $\mathscr{R}ic_f\leq (\frac{k(t)}{2h(t)}+\frac{R_g}{m+1})g$ is equivalent to $\mathscr{R}ic_f\leq (\frac{k(t)}{2h(t)}-m)g$ by the assumption that $R_g=-m(m+1)$ is preserved by the conformal Ricci flow. Note that the equation $\mathscr{R}ic_f-\frac{1}{m+1}R_g=0$ is the quasi-Einstein equation since $R_g$ is constant, which in this case can be compared with Ricci solitons.
\end{remark}
\begin{remark}
The above result (Theorem \ref{thm11}) is expected to have further applications in the setting of homogeneous $3$-manifold geometry and conformal geometry of dimension greater than $3$ in the spirit of Fischer \cite{Fi}.
\end{remark}

Lastly, we will consider some more general parabolic equation coupled with non-dynamic conformal pressure $p(t)$ 
\begin{align}\label{19}
|(\partial_t-\Delta)u|\le \bar{p}(t)(|u|+|\nabla u|)
\end{align}
along the conformal Ricci flow. The frequency $Q(t)$ for $u$ needs not to be monotone but its derivative will be suitably bounded yielding backward uniqueness of solution when $h(t)>0$ (see Theorem \ref{thm41}).

The outline of the rest part of this paper is as follows: The next section (Section \ref{sec2} is devoted to other notation and some preliminaries. The proof of main result and its applications are discussed in Section \ref{sec3}.  The last section (Section \ref{sec4}) is basically on the proof of  backward uniqueness of solution to \eqref{19}.

\section{Notation and Preliminaries}\label{sec2} 

Other notation that will be required in the sequel is presented first.  Let the dimension of the underlying manifold $M$ be a number $m+1$ not less than $3$.  Let $t$ be the abstract time parameter in the half-closed region $[0,T), T>0$. The scalar curvature, Ricci curvature and volume element of $(M,g(t))$ are respectively denoted by $R=R_{g(t)}$, $Ric=Ric(g(t))$ and $d\mu=d\mu_{g(t)}$.  In local coordinates $(x^1,\cdots,x^{m+1})$,  
$$ d\mu_{g(t)}:=\sqrt{\text{det}g_{ij}(t)}dx^1 \wedge \cdots\wedge dx^{m+1}.$$
Let $\nabla=\nabla_{g(t)}$ and $\Delta=\Delta_{g(t)}$ be the Levi-Civita connection and the Laplace-Beltrami operator with respect to $g(t)$.   Denote $|\cdot|_{g(t)} = g(t)(\cdot,\cdot)^{\frac{1}{2}}$ called $g(t)$-metric norm,  e.g., $|\nabla u|^2_{g(t)}=\langle\nabla u, \nabla u\rangle_{g(t)}$, where  $\langle\cdot, \cdot\rangle_{g(t)}$ is the inner product with respect to metric $g(t)$.

The time-dependent drifting Laplacian or $f$-Laplacian (also called weighted Laplacian) for a smooth function $f$ on $M$ is denoted by 
$$\mathscr{L}_f(\cdot):=\mathscr{L}_{g(t),f}(\cdot):= e^f\text{div}(e^{-f}\nabla(\cdot)) = \Delta(\cdot) - \langle\nabla f, \nabla(\cdot)\rangle.$$
The weighted form of the Ricci curvature tensor is the so called Bakry-\'Emery curvature tensor
$$\mathscr{R}ic_f(t):= Ric(t)+\text{Hess} f,$$
where $\text{Hess} f$ is the Hessian of function $f$. Note that having a condition of the form $\mathscr{R}ic_f=\kappa g$ is saying that $(M,g(t))$ is a Ricci soliton, which is a special solution to the Ricci flow \cite{Ha93} and very useful in the singularity analysis of the Ricci flow. These notations with and without subscript $g(t)$ are used interchangeably without resulting to any confusion. It is well know that $\mathscr{L}_f$ and $\mathscr{R}ic_f$ are related via the weighted (or drifting)  Bochner formula for an atleast $C^3$-function $h$
\begin{equation}\label{B}
\frac{1}{2} \mathscr{L}_f(|\nabla h|^2)= | \text{Hess}\ h|^{2}+\langle\nabla h, \nabla \mathscr{L}_f h\rangle+ \mathscr{R}ic_{f}(\nabla h, \nabla h).
\end{equation}

Along the flow \eqref{CRF} on a closed manifold of dimension $\ge 3$,  we let $\tau(t)=T-t$ be the backward time and define the following conjugate heat equation
on $(M^{m+1},g(t),p(t))$ as
\begin{align}
\partial_t H(t) = -\Delta_{g(t)}H(t)+(m+1)p(t)H(t)
\end{align}
with the fundamental solution
\begin{align*}
H(t)=(4\pi\tau(t))^{-\frac{m+1}{2}}e^{-f(t)}.
\end{align*}
One can then show that $f(t)$ satisfies the conjugate heat equation
\begin{align}
\partial_t f(t)= -\Delta_{g(t)}f(t)+|\nabla_{g(t)} f(t)|^2 -(m+1)\left(p(t) -\frac{1}{2\tau(t)}\right).
\end{align}
Define the weighted volume form as 
\begin{align*}
dV_{g(t)} = H(t)d\mu_{g(t)} = (4\pi\tau(t))^{-\frac{m+1}{2}}e^{-f(t)}d\mu_{g(t)}
\end{align*}
satisfying  $\int_M dV_{g(t)} =1$. Recall that $d\mu_{g(t)}$ evolves under conformal Ricci flow \eqref{CRF2} by the formula (see \cite{Fi,LLWY})
\begin{align}\label{24}
\partial_t (d\mu_{g(t)})= -(m+1)p(t)d\mu_{g(t)}.
\end{align}
Therefore one can compute that
\begin{align}\label{25}
\partial_t (dV_{g(t)}) = [H_t-(m+1)p(t)H(t)]d\mu_{g(t)} = -\frac{\Delta_{g(t)} H(t)}{H(t)}dV_{g(t)}.
\end{align}

For a smooth function $v:M\times[t_{0},t_{1}]\to \mathbb{R}$ with $v(\cdot,t),\partial_{t}v(\cdot,t)\in W_{0}^{2,2}(dV)$ for any $t\in[t_{0},t_1]\subset[0,T)$,  we define quantities $I(t)$ and $E(t)$  as follows
\begin{eqnarray}\label{26}
&&I(t)=\int_{M}v^{2}dV_{g(t)},\\\label{27}
&&E(t)=- h(t)\int_M \langle v, \mathscr{L}_fv \rangle dV_{g(t)}  =  h(t)\int_{M}|\nabla v|_{g(t)}^{2}dV_{g(t)}.
\end{eqnarray}
To this end,  referring to definitions in \eqref{26} and \eqref{27}  and reverting to \eqref{Q},  $Q(t)$ is thus written as 
\begin{eqnarray}\label{28}
 Q(t)=\frac{E(t)}{I(t)} e^{-\int_{t_{0}}^{t}\left(2\bar{p}(s)+\frac{h'(s)+k(s)}{h(s)}\right)ds},
\end{eqnarray}
where $h(t)$ and $k(t)$ are both smooth functions with respect to time-variable $t \in [0,T)$.

\section{Proof of main theorem and its applications}\label{sec3}
At first, we state some important results that will be required in the proof of the main theorem. Here and in what follows, by the virtues of \eqref{CRF2} and \eqref{lh}, the triple $(g(t),p(t),v(t))$, $t\in [0,T)$  solves the following system
\begin{align}\label{e31}
\left\{
\begin{array}{ll}
\displaystyle \partial_t g = -2(Ric(t) + (m+p(t))g(t)),  \\ 
\displaystyle (-\Delta + (m+1))p(t) = \frac{1}{m} |Ric +mg|^2,  \\
\displaystyle  \partial_t v-\Delta_{g(t)} v=  \bar{p}(t)v,  \\
\displaystyle (g(0), v(0)) = (g_0, v_0)
\end{array}
\right.
\end{align}
on closed connected manifold of dimension $m+1\ge 3$ with constant scalar curvature.  Since the conformal pressure $p(t)$ is a non-dynamical field, no initial value of $\bar{p}(t)=\max_{x\in M}p(t)$ is required.

\begin{lemma}\label{lemp} (Finiteness of $p(t)$).
Let $(g(t),p(t))$, $t\in [0,T)$ be a smooth solution to the conformal Ricci flow on a closed manifold with $R_0=-m(m+1)$ satisfying $|Ric|(x,t)\le K(t)$ for all $(x,t)\in M^{m+1}\times [0,T)$. Then the conformal pressure $p(t)$ satisfies $0\le p(t)\le K^2(t)$, for all $t\in [0,T)$, that is, $p(t)$ is finite under the flow $g(t)$.
\end{lemma} 
\proof
We know that $p(t)$ solves the elliptic equation (i.e., the second equation) in system \eqref{e31} on $M^{m+1}\times [0,T)$. So by the strong maximum principle, $p(t)\ge 0$ for all  $(x,t)\in M^{m+1}\times [0,T)$. The following conclusion can be reached: either (a) $p(t)=0$ and $Ric+mg=0$ or (b) $p(t)>0$ and $Ric+mg\neq 0$ (see \cite[Proposition 3.3]{Fi} for detail).

Suppose $(x_0,t)$ is the maxmum point, that is $p(x_0,t):=\max_{x\in M}p(x,t)$, we have $\nabla p(x_0,t)=0$ and $\Delta p(x_0,t)\le 0$. Hence
\begin{align*}
m(m+1)p(x_0,t)\le |Ric+mg|^2 (x_0,t)= |Ric|^2(x_0,t) -m^2(m+1) \le K^2(t).
\end{align*}
Hence,  $0\le p(x,t)\le K^2(t)$ for all $(x,t)\in M^{m+1}\times [0,T)$. 

\qed

\begin{lemma}\label{lem31}
Let $(g(t),p(t),v(t))$, $t\in [0,T)$ solves the system \eqref{e31}. The following identities hold:
\begin{align}\label{e32}
\partial_t \left(|\nabla v|^2_{g(t)}\right) = 2 (Ric+mg(t))(\nabla v, \nabla v) +2(p(t)+ \bar{p}(t))|\nabla v|^2_{g(t)} +2\langle\nabla v, \nabla\Delta v\rangle,
\end{align}
\begin{align}\label{e33}
\left(\partial_t - \Delta_{g(t)}\right) |\nabla v|^2_{g(t)} = 2(p(t)+\bar{p}(t)) |\nabla v|^2_{g(t)} +2mg(t)(\nabla v, \nabla v)- 2|\text{Hess}\ v|^2_{g(t)},
\end{align}
where $\bar{p}(t):=\max_{x\in M}p(t)$.
\end{lemma}

\proof
Following the standard computation under geometric flow we have 
\begin{align*}
\partial_t (|\nabla v|^2) = -[\partial_t g] (\nabla v, \nabla v) + 2\langle\nabla v, \nabla \partial_t v\rangle.
\end{align*}
Reverting to system (\ref{e31}), we have values for the quantities $-[\partial_t g]$ and $\partial_t v$, which when substituted into the last expression together with the fact that $\nabla \bar{p}(t)=0$ yields \eqref{e32}.  
Combining (\ref{e31}) with the classical Bochner formula proves (\ref{e33}).

\qed
 
\begin{lemma}\label{lem32}
For all $u,v \in W^{1,2}_0(dV_{g(t)})$, the drifting Laplacian $\mathscr{L}_{f(t)}$ satisfies integration by parts formula, i.e.,
\begin{align*}
\int_Mu \mathscr{L}_{f(t)} v dV_{g(t)} = - \int_M\langle \nabla u, \nabla v \rangle_{g(t)} dV_{g(t)},
\end{align*}
and it is self-adjoint with respect to the weighted measure $dV_{g(t)} $, i.e.,
\begin{align*}
\int_Mu \mathscr{L}_{f(t)} v dV_{g(t)} = \int_M (\mathscr{L}_{f(t)} u)v dV_{g(t)}.
\end{align*}
\end{lemma} 

\proof
Recall that $\mathscr{L}_{f(t)} v := e^{f(t)}\text{div}(e^{-f(t)}\nabla v)$ and $dV_{g(t)}:= (4\pi\tau(t))^{-\frac{m+1}{2}} e^{-f(t)}d\mu_{g(t)}$.
Direct computation using classical integration by parts gives
\begin{align*}
\int_Mu\mathscr{L}_fv\ dV &= (4\pi\tau)^{-\frac{m+1}{2}}\int_Mue^f\text{div}(e^{-f}\nabla v)e^{-f}d\mu\\
&= - (4\pi\tau)^{-\frac{m+1}{2}}\int_M  e^{-f} \langle \nabla u, \nabla v \rangle d\mu := - \int_M   \langle \nabla u, \nabla v \rangle dV \\
& = (4\pi\tau)^{-\frac{m+1}{2}}\int_M \text{div}(e^{-f}\nabla u)v d\mu\\
& = (4\pi\tau)^{-\frac{m+1}{2}}\int_M e^f \text{div}(e^{-f}\nabla u)v e^{-f}d\mu := \int_M (\mathscr{L}_fu)v\ dV.
\end{align*}
That is, 
\begin{align*}
\int_Mu\mathscr{L}_fv\ dV  = - \int_M   \langle \nabla u, \nabla v \rangle dV = \int_M (\mathscr{L}_fu)v\ dV.
\end{align*}
\qed 
 
\begin{lemma}\label{lem33} (Drifting Reilly formula).
For any $v\in W^{2,2}_0(dV_{g(t)})$, then
\begin{align}\label{e34}
\int_M |\text{Hess}\ v|^2_{g(t)} = \int_M \left[ (\mathscr{L}_{f(t)}v)^2 -\mathscr{R} ic_f(\nabla v, \nabla v)\right]dV_{g(t)}.
\end{align}
\end{lemma} 
\proof
Applying the drifting Bochner formula \eqref{B}, integration by parts and self-adjoint properties of $\mathscr{L}_f$ with respect to the weighted measure $dV$, we have
\begin{align*}
0=\frac{1}{2}\int_M \mathscr{L}_f(|\nabla v|^2)\ dV = \int_M\left[ |\text{Hess}\ v|^2_{g(t)} -  (\mathscr{L}_fv)^2  + \mathscr{R}ic_f(\nabla v, \nabla v)\right]dV
\end{align*}
which is the required formula.
\qed 
 
 \subsection*{Proof of Theorem \ref{thm11}}
We only give the proof of Part (i) since the proof of Part (ii) is similar. The procedure is to compute derivatives of $I(t)$ and $E(t)$, and then use Cauchy-Schwarz inequality to show that $Q'(t)$ is nonpositive.
 
 \proof
 Suppose $v$ solves the heat equation in the system \eqref{e31}, i.e., $v_t=\Delta v+\bar{p}(t)v$.  In consideration of evolution of the weighted measure in \eqref{25}, integration by parts formula and the identity $\Delta v^2=2v\Delta v+ 2|\nabla v|^2$, we first compute the derivative of $I(t)$ (recall that $I(t)$ is defined in \eqref{26}) as follows
\begin{align*}
I'(t) & = \int_M\left(2vv_t - v^2 \frac{\Delta H}{H}\right)dV\\
&= \int_M\left(2v(v_t-\Delta v)-2|\nabla v|^2\right)dV \\
& =2\bar{p}(t)\int_Mv^2dv - 2\int_M |\nabla v|^2 dV = 2\bar{p}I(t) - \frac{2}{h(t)}E(t).
\end{align*} 
Similarly, we can compute derivative of the energy $E(t)$ (recall that $E(t)$ is defined in \eqref{27}) as follows
\begin{align}\label{e35}
E'(t)=h'(t)\int_M|\nabla v|^2 dV + h(t)\int_M(\partial_t -\Delta)|\nabla v|^2 dV.
\end{align}
Applying \eqref{e33} of Lemma \ref{lem31} into \eqref{e35}, we obtain 
\begin{align*}
E'(t)\le (h'(t)+4h(t)\bar{p}(t))\int_M|\nabla v|^2dV + 2 h(t)\int_M \left[mg(t)(\nabla v, \nabla v) - |\text{Hess}\ v|^2\right]dV
\end{align*} 
 due to the condition $h(t)>0$ and the fact that $p(t)\le \bar{p}(t)$.  
 Combining the last equation with drifting Reilly formula of Lemma \ref{lem33} we get
\begin{align*}
E'(t) \le (h'(t)&+4h(t)\bar{p}(t))  \int_M|\nabla v|^2dV \\
&  -2 h(t)\int_M \left[ |\mathscr{L}_fv|^2 - (\mathscr{R}ic_f + mg(t))(\nabla v, \nabla v)\right]dV,
\end{align*}  
which leads to the following inequality
\begin{align*}
E'(t) \le \left[4\bar{p}(t)+\frac{k(t)+h'(t)}{h(t)}\right] E(t) -2 h(t)\int_M|\mathscr{L}_f v|^2 dV
\end{align*}
by invoking the Bakry-\'Emery curvature bound $\mathscr{R}ic_f\leq (\frac{k(t)}{2h(t)}+\frac{R_g}{m+1})g$.
 Proceeding from here, we can express the frequency function $Q(t)$ in terms of rescaled quantities $\widetilde{I}(t)$ and $\widetilde{E}(t)$ as follows
\begin{align*}
Q(t) = \frac{\widetilde{E}(t)}{\widetilde{I}(t)}
\end{align*}
 such that
  $$\widetilde{I}(t) = I(t)e^{-\int_{t_0}^t 2\bar{p}(s)ds}$$
  and
 \begin{align*}
 \widetilde{E}(t) = E(t)e^{-\int_{t_0}^t \left(4\bar{p}(s)+\frac{k(s)+h'(s)}{h(s)}\ \right)ds}
 \end{align*}
 with their respective derivatives computed as follows
$$\widetilde{I}'(t) = -\frac{2}{h(t)}E(t)e^{-\int_{t_0}^t 2\bar{p}(s)ds}$$
  and
 \begin{align*}
 \widetilde{E}'(t) \le  -2h(t) e^{-\int_{t_0}^t \left(4\bar{p}(s)+\frac{k(s)+h'(s)}{h(s)} \right)ds} \int_M|\mathscr{L}_f v|^2 dV.
 \end{align*}
Now using the bound for $ \widetilde{E}(t) $ will allow the derivative of the frequency functional $Q(t)$ to be bounded. Suppose $h(t)>0$ and denote by $\Gamma(t)$ the following integral
$$\Gamma(t):= \int_{t_0}^t \left(6\bar{p}(s)+\frac{k(s)+h'(s)}{h(s)} \right)ds.$$
 \begin{align*}
 \widetilde{I}^2(t) Q'(t) & =  \widetilde{I}(t) \widetilde{E}'(t) - \widetilde{I}'(t) \widetilde{E}(t) \\
 & \le e^{-\Gamma(t)} \left[-2h(t)I(t)\int_M(\mathscr{L}_f v)^2 dV +  \frac{2}{h(t)}E^2(t) \right]\\
 & =  -2h(t) e^{-\Gamma(t)} \left[I(t)\int_M(\mathscr{L}_f v)^2 dV -\left( \frac{1}{h(t)}E(t)\right)^2 \right]\\
& = -2h(t) e^{-\Gamma(t)} \left[\left(\int_M v^2dV\right)\left(\int_M(\mathscr{L}_f v)^2 dV\right) -\left(\int_M|\nabla v|^2dV\right)^2 \right]\\
&\le 0.
 \end{align*}
The first inequality is due to the bound for  $\widetilde{E}(t)$ whilst the last inequality is due to integration by parts formula and Cauchy-Schwarz inequality. That is,
 \begin{align*}
 \int_M|\nabla v|^2dV &= - \int_M \langle v, \mathscr{L}_fv\rangle dV \le \int_M |v||\mathscr{L}_f v| dV\\
 &\le \left(\int_M |v|^2dV\right)^{\frac{1}{2}}\left(\int_M |\mathscr{L}_f v|^2 dV\right)^{\frac{1}{2}}
 \end{align*}
 which implies
 \begin{align}\label{e37}
 \left(\int_M |v|^2dV\right)\left(\int_M |\mathscr{L}_f v|^2 dV\right)- \left(\int_M|\nabla v|^2dV\right)^2\ge 0.
 \end{align}
 Hence, $Q(t)$ is a nonincreasing function along the conformal Ricci flow if $h(t)>0$. The proof for the case of $h(t)<0$ follows suit. 
 
 Moreover, suppose $Q'(t)=0$, the equality in Cauchy schwarz inequality \eqref{e37} implies $-\mathscr{L}_fv=c(t)v$, where $c(t) = \frac{1}{h(t)}Q(t) e^{\int_{t_0}^t \left(2\bar{p}(s)+'\frac{h'(s)+k(s)}{h(s} \right) ds}$.
 
 \qed
 
 \subsection*{Proof of Corollary \ref{cor2}}
 Recall that $I'(t)= 2\bar{p}(t)I(t) - \frac{2}{h(t)}E(t)$. Therefore
 \begin{align*}
 \frac{d}{dt}\left(\log I(t)\right) = \frac{I'(t)}{I(t)} &= 2\bar{p}(t) - \frac{2}{h(t)}\frac{E(t)}{I(t)}\\
 & =2\bar{p}(t)-  \frac{2}{h(t)}Q(t)e^{\int_{t_0}^t \left(2\bar{p}(s)+\frac{h'(s)+k(s)}{h(s)} \right) ds}.
 \end{align*}
 Integrating both sides of the last equation on $[a,b]\subset [t_0,t_1)$ and using the monotonicity of $Q$ in Theorem \ref{thm11}, we have
 \begin{align*}
 \log I(b)-\log I(a) & \ge  2\int_a^b \bar{p}(t)dt - 2 \int_a^b \frac{Q(t)}{h(t)} e^{\int_{t_0}^t \left(2\bar{p}(s)+\frac{h'(s)+k(s)}{h(s)} \right) ds}dt \\
 & \ge   2\int_a^b \bar{p}(t)dt - 2 Q(a)\int_a^b \frac{1}{h(t)} e^{\int_{t_0}^t \left(2\bar{p}(s)+\frac{h'(s)+k(s)}{h(s)} \right) ds}dt.
\end{align*}  
 Exponentiating yields (since $p(t)$ and $h(t)$ are finite)
 \begin{align*}
 \frac{I(b)}{I(a)} \ge \exp \left\{2\int_a^b \bar{p}(t)dt - 2 Q(a)\int_a^b \frac{1}{h(t)} e^{\int_{t_0}^t \left(2\bar{p}(s)+\frac{h'(s)+k(s)}{h(s)} \right) ds}dt \right\}.
 \end{align*}
 Therefore, if $v(\cdot,b)=0$, then $I(b)=0$ and the last inequality implies $I(a)=0$. Since $a$ is arbitrary we conclude that $I(t)=0$ for any $t\in [a,b]\subset(0,T)$. Hence $v(\cdot,t)= 0$ for any $t\in [a,b]\subset(0,T)$.
 
 \qed
 
 \subsection*{Proof of Corollary \ref{cor3}} 
 Given the time interval $[t_0,t]\subset (0,T)$. Let $(g(t),p(t),v(t))$ $t\in[0,T)$ solves \eqref{e31}, that is, $v(t)$ solves $v_t=\Delta v+\bar{p}(t)v$. By the hypothesis, $v(\cdot,t_0)$ is the eigenfunction of $-\mathscr{L}_{f(t_0)}$ corresponding to the first eigenvalue $\lambda(t_0):=\lambda(g(t_0))$, thus; 
$ -\mathscr{L}_{f(t_0)} v(\cdot,t_0) = \lambda(t_0) v(\cdot,t_0)$. Then we  have the frequency functional $Q(t)$ for $v(t)$ by \eqref{28} as follows
\begin{align*}
Q(t):= \frac{- h(t)\int_M \langle v, \mathscr{L}_fv \rangle dV_{g(t)} }{\int_M v^2dV} e^{-\int_{t_{0}}^{t}\left(2\bar{p}(s)+\frac{h'(s)+k(s)}{h(s)}\right)ds}.
\end{align*}
 Based on Theorem \ref{thm11}, we know that $Q(t)$ is nonincreasing for $h(t)>0$ and nondecreasing for $h(t)<0$. To this end, we  refer to the definition of $\lambda(t):=\lambda(g(t))$ (see the statement of the corollary) and have for any $h(t)>0$ and any $t\in [t_0,t_1]\subset (0,T)$ that
 \begin{align*}
 h(t_0)\lambda(t_0) =Q(t_0)\ge Q(t) \ge h(t)\lambda(t)e^{-\int_{t_{0}}^{t}\left(2\bar{p}(s)+\frac{h'(s)+k(s)}{h(s)}\right)ds}.
 \end{align*}
As a consequence of this, it is clear that $h(t)\lambda(t)$ is monotone decreasing for $h(t)>0$.  The second part of the corollary can be proved in a similar way.

\qed

 \section{General parabolic equations}\label{sec4}
In this section we consider some more general parabolic equation coupled with non-dynamic conformal pressure $p(t)$ along the conformal Ricci flow. Here the parabolic  frequency $Q(t)$ in relation to solution $v$ of 
\begin{align}\label{e41}
|(\partial_t-\Delta)v|\le \bar{p}(t)(|v|+|\nabla v|),
\end{align}
where $\disp \bar{p}(t)=\max_{x\in M}p(t)$,  along the conformal Ricci flow needs to be monotone but its derivative can be suitably bounded to imply  backward uniqueness of solution when $h(t)>0$. The main theorem of this section is the following:
 
 \begin{theorem}\label{thm41}
 Let $v:M^{m+1}\times [t_0,t_1]\to \mathbb{R}$, $t_0<t_1$ satisfy \eqref{e41} along the conformal Ricci flow $(M^{m+1},g(t),p(t))$, $t\in [t_0,t_1]$ with $\mathscr{R}ic_f\leq (\frac{k(t)}{2h(t)}+\frac{R_g}{m+1})g$ and $R_g=-m(m+1)$. If $v(\cdot,t_1)=0$, then $v(\cdot,t)\equiv 0$ for all $t\in [t_0,t_1]$.
 \end{theorem}
 
As in Section \ref{sec2},  define $\tau(t)=T-t$ to be the backward time and the conjugate fundamental solution at some point $(x,t)$ to be
\begin{align*}
H(t)=(4\pi\tau(t))^{-\frac{m+1}{2}}e^{-f(t)}
\end{align*}
on $(M^{m+1},g(t),p(t))$  so that the weighted measure remains $dV =Hd\mu$.

For a smooth function $v:M\times[t_{0},t_{1}]\to \mathbb{R}$ with $v(\cdot,t),\partial_{t}v(\cdot,t)\in W_{0}^{2,2}(dV)$ and for any $t\in[t_{0},t]\subset[0,T)$,  we define the quantities $I(t)$, $E(t)$ and the parabolic functional $Q(t)$ as in \eqref{26}--\eqref{28}.    
 
\begin{proposition}\label{prop42}
 Let $v:M^{m+1}\times [t_0,t_1]\to \mathbb{R}$ satisfy \eqref{e41} along the conformal Ricci flow $(M^{m+1},g(t),p(t))$, $t\in [t_0,t_1]$ with $\mathscr{R}ic_f\leq (\frac{k(t)}{2h(t)}+\frac{R_g}{m+1})g$ and $R_g=-m(m+1)$.  Then
 \begin{align}\label{e42}
 (\log I(t))' \ge -3\bar{p}(t)-\frac{\bar{p}(t)+2}{h(t)}Q(t) e^{\int_{t_0}^t \left(2\bar{p}(s)+\frac{h'(s)+k(s)}{h(s)} \right) ds},
 \end{align}
  \begin{align}\label{e43}
 Q'(t)\le \bar{p}^2(t)\left[Q+h(t_0)\right]
 \end{align}
 and
  \begin{align}\label{e44}
\bar{p}^2(t)\ge [\log(Q+h(t_0)]'.
 \end{align}
\end{proposition} 

\proof
By direct differentiation and integration by parts formula we have
\begin{align*}
I'(t) & = \int_M\left(2vv_t - v^2 \frac{\Delta H}{H}\right)dV= \int_M\left(2v(v_t-\Delta v)-2|\nabla v|^2\right)dV \\
& \ge - 2\bar{p}(t)\int_M|v|(|v|+ |\nabla v|) dV -2\int_M |\nabla v|^2 dV\\
&=-2\bar{p}I(t) - 2\bar{p}(t)\int_M|v||\nabla v|dV- \frac{2}{h(t)}E(t)\\
&\ge -3\bar{p}(t)I(t) - \bar{p}(t) \int_M |\nabla v|^2dV- \frac{2}{h(t)}E(t)\\
&= -3\bar{p}(t)I(t) -  \frac{\bar{p}(t)+2}{h(t)}E(t) \\
& = -3\bar{p}(t)I(t) -  \frac{\bar{p}(t)+2}{h(t)}I(t)Q(t) e^{\int_{t_0}^t \left(2\bar{p}(s)+\frac{h'(s)+k(s)}{h(s)} \right) ds},
\end{align*} 
where we have used \eqref{e41} to get the first inequality, and an elementary inequality of the form $2rs\le r^2+s^2$ to get the second inequality. This prove \eqref{e42}.

Applying integration by parts formula again we have
\begin{align*}\label{e45}
I'(t) & = \int_M\left(2vv_t - v^2 \frac{\Delta H}{H}\right)dV\\
& = \int_M\left( 2vv_t -2v\Delta v-2|\nabla v|^2\right)dV \\
&=2\int_Mv\left(\partial_t-\Delta + \mathscr{L}_f\right)vdV\\
& = 2 \int_Mv\left[\mathscr{L}_f+\frac{1}{2}(\partial_t- \Delta) \right]vdV + \int_M v(\partial_t- \Delta)vdV.
\end{align*} 
Letting
\begin{align*}
\widetilde{I}(t) = I(t)e^{-\int_{t_0}^t2\bar{p}(s)ds}.
\end{align*}
Then
\begin{align}
\widetilde{I}'(t) & = \Bigg\{2 \int_Mv\left[\mathscr{L}_f+\frac{1}{2}(\partial_t- \Delta) \right]vdV + \int_M v(\partial_t- \Delta)vdV\\
& \hspace{2cm}  -2\bar{p}(t)I(t) \Bigg\}e^{-\int_{t_0}^t 2\bar{p}(s)ds}. \nonumber
\end{align}

Note that we can  also write $E(t)$ as follows
\begin{align}\label{e46}
E(t) & = -h(t)\int_Mv\mathscr{L}_fvdV \nonumber\\
& = -h(t)\int_Mv\left[\mathscr{L}_f+\frac{1}{2}(\partial_t- \Delta) \right]vdV + \frac{h(t)}{2}\int_M v(\partial_t- \Delta)vdV.
\end{align}
Along the conformal Ricci flow, and with the aid of the classical Bochner formula, we have 
\begin{align}\label{e47}
(\partial_t-\Delta)|\nabla v|^2 & = 2mg(t)(\nabla v,\nabla v)+2p(t)|\nabla v|^2 \nonumber\\
& \hspace{1cm} -2|\text{Hess}\ v|^2 + 2\langle\nabla v, \nabla (\partial_t- \Delta)v\rangle.
\end{align}
Computing the derivative of $E(t)$ using \eqref{e47} we have
\begin{align*}
E'(t) & = h'(t)\int_M|\nabla v|^2dV +h(t)\int_M(\partial_t-\Delta)|\nabla v|^2dV\\
&= h(t)\left[\frac{h'(t)}{h(t)} \int_M|\nabla v|^2dV + \int_M(\partial_t-\Delta)|\nabla v|^2dV\right]\\
& = h(t)\int_M \Big[ \frac{h'(t)+2h(t)p(t)}{h(t)}|\nabla v|^2 +2mg(t)(\nabla v,\nabla v) \\
& \hspace{2cm} -2|\text{Hess}\ v|^2 + 2\langle\nabla v, \nabla (\partial_t- \Delta)v\rangle \Big]dV.
\end{align*}
Applying integration by parts formula and drifting Reilly formula (Lemma \ref{lem33} gives
\begin{align*}
E'(t) &  = h(t)\int_M \Big[ \frac{h'(t)+2h(t)p(t)}{h(t)}|\nabla v|^2 -2 |\mathscr{L}_fv|^2 +2(\mathscr{R}ic_f +mg)(\nabla v,\nabla v) \\
& \hspace{2cm}  - 2(\mathscr{L}_fv) (\partial_t- \Delta)v\Big]dV.
\end{align*}
Applying the Bakry-\'Emery Ricci curvature bound  $\mathscr{R}ic_f\leq (\frac{k(t)}{2h(t)}+\frac{R_g}{m+1})g$ with $R_g=-m(m+1)$ gives
\begin{align*}
E'(t) & \le -2h(t)\int_M \left[  |\mathscr{L}_fv|^2 + (\mathscr{L}_fv) (\partial_t- \Delta)v \right]dV \\
&\hspace{1cm} + (k(t)+h'(t)+2h(t)\bar{p}(t))\int_M|\nabla v|^2dV\\
& = -2h(t)\int_M \left[ \left| \left(\mathscr{L}_f+\frac{1}{2}\Big(\partial_t-\Delta\Big)\right) v \right|^2 -\frac{1}{4}\Big|\Big(\partial_t-\Delta\Big)v \Big|^2 \right] dV\\
&\hspace{1cm} + \left[ 2\bar{p}(t)+\frac{h'(t)+k(t)}{h(t)}\right]E(t),
\end{align*}
where we have used the completing the square method. Also letting
\begin{align}\label{e48a}
\widetilde{E}(t)= E(t)e^{-\int_{t_0}^t\left(4\bar{p}(s)+\frac{h'(s)+k(s)}{h(s)}\right)ds}.
\end{align}
Then
\begin{align*}
\widetilde{E}'(t) = E'(t) & e^{-\int_{t_0}^t\left(4\bar{p}(s)+\frac{h'(s)+k(s)}{h(s)}\right)ds} \\
& - \left( 4\bar{p}(t)+\frac{h'(t)+k(t)}{h(t)} \right)E(t)e^{-\int_{t_0}^t\left(4\bar{p}(s)+\frac{h'(s)+k(s)}{h(s)}\right)ds}.
\end{align*}
Using the above bound for $E'(t)$ gives (since $p(t)\le \bar{p}(t)$)
\begin{align}\label{e49}
\widetilde{E}'(t) & \le  \Bigg\{-2h(t)\int_M \left[ \left| \left(\mathscr{L}_f+\frac{1}{2}\Big(\partial_t-\Delta\Big)\right) v \right|^2 -\frac{1}{4}\Big|\Big(\partial_t-\Delta\Big)v \Big|^2 \right] dV\\
&\hspace{1cm} -2\bar{p}(t)E(t)\Bigg\}e^{-\int_{t_0}^t\left(4\bar{p}(s)+\frac{h'(s)+k(s)}{h(s)}\right)ds}.\nonumber
\end{align}
Combining \eqref{e45}, \eqref{e46} and \eqref{e48a} we get 
\begin{align}\label{e410a}
\widetilde{I}'(t)\widetilde{E}(t) &=  \Bigg\{-2h(t)\left( \int_Mv\left[\mathscr{L}_f+\frac{1}{2}(\partial_t- \Delta) \right]vdV\right)^2 \nonumber \\
&+ \frac{h(t)}{2}\left(\int_M v(\partial_t- \Delta)vdV\right)^2  -2\bar{p}(t)I(t)E(t) \Bigg\}e^{-\int_{t_0}^t\left(6\bar{p}(s)+\frac{h'(s)+k(s)}{h(s)}\right)ds}.
\end{align}
Now, using the bound for $\widetilde{E}'(t)$ (i.e., \eqref{e49}), \eqref{e410a}, Cauchy-Schwarz inequality and elementary inequality of the form $(r+s)^2\le2(r^2+s^2)$,  with $\Gamma_2(t)$ denoting
$$\Gamma_2(t):=\int_{t_0}^t\left(6\bar{p}(s)+\frac{h'(s)+k(s)}{h(s)}\right)ds$$
we obtain
\begin{align*}
& \widetilde{I}^2(t)Q'(t) = \widetilde{I}(t)\widetilde{E}'(t)-\widetilde{I}'(t)\widetilde{E}(t)\\
& \le e^{-\Gamma_2(t)} \Bigg\{-2h(t)I(t)\int_M \left[ \left| \left(\mathscr{L}_f+\frac{1}{2}\Big(\partial_t-\Delta\Big)\right) v \right|^2 -\frac{1}{4}\Big|\Big(\partial_t-\Delta\Big)v \Big|^2 \right] dV\\
& \hspace{1cm} +2h(t) \left( \int_Mv\left[\mathscr{L}_f+\frac{1}{2}(\partial_t- \Delta) \right]vdV\right)^2 -  \frac{h(t)}{2}\left(\int_M v(\partial_t- \Delta)vdV\right)^2\Bigg\}\\
&=-2h(t)e^{-\Gamma_2(t)} \Bigg\{I(t)\left( \int_M \left| \left(\mathscr{L}_f+\frac{1}{2}\Big(\partial_t-\Delta\Big)\right) v \right|^2 \right) -  \left( \int_Mv\left[\mathscr{L}_f+\frac{1}{2}(\partial_t- \Delta) \right]vdV\right)^2\\
& \hspace{1cm} + \frac{1}{4}\left(\int_M v(\partial_t- \Delta)vdV\right)^2 -\frac{1}{4}I(t) \int_M \Big|\Big(\partial_t-\Delta\Big)v \Big|^2dV\Bigg\}\\
&\le \ \  \frac{h(t)}{2}I(t) \left( \int_M \Big|\Big(\partial_t-\Delta\Big)v \Big|^2dV \right)e^{-\Gamma_2(t)} \\
&\le \ \  \frac{h(t)}{2}\bar{p}^2(t)I(t) \left( \int_M  \Big(|v|+|\nabla v|\Big)^2dV \right)e^{-\Gamma_2(t)} \\
&\le \ \  \bar{p}^2(t)I(t)  \Big(h(t)I(t)+E(t)\Big)e^{-\Gamma_2(t)}.
\end{align*}
Using the condition $h(t)>0$ we have 
\begin{align*}
Q'(t) &\le \bar{p}^2(t)\left(h(t)+\frac{E(t)}{I(t)}\right)e^{-\int_{t_0}^t\left(2\bar{p}(s)+\frac{h'(s)+k(s)}{h(s)}\right)ds}\\
& = \bar{p}^2(t)\left(Q(t) + h(t) e^{-\int_{t_0}^t\left(2\bar{p}(s)+\frac{h'(s)+k(s)}{h(s)}\right)ds}\right)\\
&\le \bar{p}^2(t)\Big(Q(t)+h(t_0)\Big)
\end{align*}
which proves \eqref{e43} and \eqref{e44} of the proposition.

\qed

\begin{corollary}
Let $v:M^{m+1}\times [t_0,t_1]\to \mathbb{R}$ satisfy \eqref{e41} along the conformal Ricci flow $(M^{m+1},g(t),p(t))$, $t\in [t_0,t_1]\subset [0,T)$. Then
\begin{align}\label{e411}
I(t_1)\ge & I(t_0)\exp \Bigg\{-3(t_1-t_0)\sup_{t\in [t_0,t_1]}\max_{x\in M} p(t)- \Bigg[\left(2+\sup_{t\in [t_0,t_1]}\max_{x\in M} p(t)\right)\nonumber \\
&\times \Big(Q(t_0)+h(t_0)\Big) \left(e^{\int_{t_0}^{t_1}\bar{p}^2(t)dt }\right)\Bigg] \int_{t_0}^{t_1}\frac{1}{h(t)}e^{\int_{t_0}^{t_1} 2\bar{p}(s)+\frac{h'(s)+k(s)}{h(s)} ds} dt\Bigg\}.
\end{align}
Moreover,  if $v(\cdot,t_1)=0$, then $v(\cdot,t)\equiv 0$ for all $t\in [t_0,t_1]$.
\end{corollary}

This corollary is indeed a restatement of Theorem \ref{thm41}. 

\proof
Integrating \eqref{e42} on the interval $[t_0,t_1]$ gives 
\begin{align}\label{e412}
 \log I(t_1)  - & \log I(t_0) \nonumber \\
& \ge  -3\int_{t_0}^{t_1}\bar{p}(t)dt  - \int_{t_0}^{t_1} \frac{\bar{p}(t)+2}{h(t)}Q(t) e^{\int_{t_0}^t \left(2\bar{p}(s)+\frac{h'(s)+k(s)}{h(s} \right) ds} dt \nonumber\\
&\ge -3(t_1-t_0)\sup_{t\in [t_0,t_1]}\bar{p}(t)\\
& \hspace{1cm} - \left(2+\sup_{t\in [t_0,t_1]}\bar{p}(t)\right)\int_{t_0}^{t_1} \frac{Q(t)}{h(t)} e^{\int_{t_0}^t \left(2\bar{p}(s)+'\frac{h'(s)+k(s)}{h(s} \right) ds} dt. \nonumber
\end{align}
Integrating \eqref{e44} we get 
\begin{align*}
\log[Q(t)+h(t_0)] \le \log[Q(t_0)+h(t_0)] + \int_{t_0}^{t_1}\bar{p}^2(t)dt.
\end{align*}
Therefore, $Q(t)$ is bounded by
\begin{align}\label{e413}
Q(t)\le (Q(t_0)+h(t_0))e^{\int_{t_0}^{t_1}\bar{p}^2(t)dt} -h(t_0).
\end{align}
Inserting the bound \eqref{e413} into \eqref{e412},  then \eqref{e411} is obtained by exponentiation.

The conclusion in the second part of the corollary follows immediately from the first part. Note that the fact that $p(t)\ge 0$ is finite and the assumption that $h(t)>0$ (uniformly) for all $t\in [t_0,t_1]\subset [0,T)$ makes the integral appearing in \eqref{e411} finite.

\qed

\section*{Conflict of interest}
The authors declare that there is no conflict of interest.


\begin{thebibliography}{99}

\bibitem{Am}  
F. J. Almgren, 
{\it Dirichlet's problem for multiple valued functions and the regularity of mass minimizing integral currents}.  Minimal submanifolds and geodesics (M. Obata Ed.), 1-6, North Holland, Amsterdam, 1976.

\bibitem{Au} 
T. Aubin,   
{\it Some nonlinear problems in Riemannian geometry}, 
Springer Monographs in Mathematics, Springer (1998).

\bibitem{AA} 
S.  Azami and A. Abolarinwa,
 {\it Parabolic Frequency on Ricci-Bourguignon flow  and Yamabe flow}. 
 Submitted (2023).
 
\bibitem{AA1} 
S.  Azami and A. Abolarinwa,
 {\it  Yamabe constant evolution along the Ricci-Bourguignon flow}.  Arab. J. Math. , 11 (2022),  459--467.
 
\bibitem{BJ} 
J. Baldauf and D. Kim 
{\it Parabolic frequency on Ricci flows},  
Int. Math. Res. Notice,  2022; rnac128, https://doi.org/10.1093/imrn/rnac128.

\bibitem{BL} 
J. Baldauf and T.-K. Lee
{\it Parabolic frequency for the mean curvature flow},  
https://doi.org/10.48550/arXiv.2210.14286

\bibitem{BB} 
N. Basu and A. Bhattacharyya,
{\it Evolution of $\mathcal{F}$-functional and $\omega$-entropy functional for the conformal Ricci flow},  
Acta Univ. Sapientae Math. , 6(2) (2014),  209-216.

\bibitem{Bel} 
T. Bell,
{\it Uniqueness of conformal Ricci flow using energy method},  
Pacific J. Math. 286 (2) (2017),  277-290.

  \bibitem{CM1} 
T. H. Colding and W. P. Minicozzi II, 
{\it Harmonic functions with polynomial growth},
 J. Diff.  Geom., 46(1),  (1997), 1-77.

\bibitem{CM2} 
T. H. Colding and W. P. Minicozzi II, 
{\it Parabolic frequency on manifolds},
 Int.  Math. Res.  Notice,  2022 (15) (2022), 11878-11890.
 
 
\bibitem{Fi} A. E. Fischer, {\it An introduction to conformal Ricci flow}, Classical Quantum Gravity 21(3)(2004), 171-218.

\bibitem{GL1} 
N. Garofalo and F. H. Lin, 
{\it Monotonicity properties of variational integrals, $A_p$ weights and unique continuation},
 Indiana Univ. Math. J., 35(2) (1986), 245-268.
 
\bibitem{GL2} 
N. Garofalo and F. H. Lin, 
{\it Unique continuation for elliptic operators: a geometric-variational approach}, Comm. Pure Appl. Math., 40(3) (1987), 347-366.



\bibitem{Ha82} 
R. S. Hamilton, 
{\it Three-manifolds with positive Ricci curvature},  Journal of Differential Geometry. 17 (2), (1982), 255-306

\bibitem{Ha93} 
R.S. Hamilton, 
{\it Monotonicity formulas for parabolic flows on manifolds},  Communications in Analysis and Geometry,  1(1),  (1993),  127-137.

\bibitem{HL}
Q. Han and F. H.  Lin, 
{\it Nodal sets of solutions of parabolic equations: II},  
Comm. on Pure Appl. Math. , 47(9), (1994), 1219-1238.


\bibitem{LASA}
Y. Li, A. Abolarinwa, S. Azami and A. Ali,
{\it Yamabe constant evolution and monotonicity along the conformal Ricci flow},
AIMS Math., 7(7), (2022), 12077-12090.

\bibitem{LLX} 
C. Li, Y. Li and K. Xu,  
{\it Parabolic frequency monotonicity on Ricci-flow and Ricci-harmonic flow with bounded curvature}, 
  arxiv.org/abs/2205.07702v1

\bibitem{LLWY}  
F Li,  P. Lu, J. Wang and Y. Zheng 
{\it Monotonicity of functionals along conformal Ricci flow},
Proc. Amer. Math. Soc.,  148(9), (2020), 4007-4014.



\bibitem{Li} 
F. H. Lin, 
{\it Nodal sets of solutions of elliptic and parabolic equations},
 Comm. Pure Appl. Math.,  44(3) (1991), 287-308.
 
\bibitem{Lo}
A. Logunov, 
{\it Nodal sets of Laplace eigenfunctions: polynomial upper estimates of
the Hausdorff measure.},
 Annals of Mathematics,  (2018), 221-239.
  
 

\bibitem{LQZ} P. Lu, J Qing and Y. Zheng, 
{\it A note on conformal Ricci flow}, 
Pacific J. Math. 268 (2014), 413-434.

\bibitem{LQZ1} P. Lu, J Qing and Y. Zheng, 
{\it Conformal Ricci flow on asymptotically hyperbolic manifolds}, 
Sci.  China Math.  62(1), (2019),  157-170.
 

\bibitem{SZ} 
X. Sun and A. Zhu, 
{\it Backward uniqueness for the conformal Ricci flow}, 
Diff. Geom.  Appl.  56 (2018),  110-119.  ,

\bibitem{N}
 L. Ni,
 {\it Parabolic frequency monotonicity and a theorem of Hardy-P\'olya-Szeg\"o},  Analysis, complex geometry, and mathematical physics: in honor of Duong H. Phong,  203-210.  Contemp.  Math., 644, Amer. Math. Soc., Providence,  RI, 2015.


\bibitem{P} 
C. C. Poon, 
{\it Unique continuation for parabolic equations},  
Comm. Partial Differ.  Equ., 21(3-4) (1996),  521-539.



\bibitem{XL} 
X. Li and K.  Wang,  
{\it Parabolic frequency monotonicity on compact manifolds},
Calc. Var.  58 (2019), no. 6, Paper No. 189, 18 pp.


\end{thebibliography}
\end{document}